\documentclass{article}

\usepackage{upref}
\usepackage{amsfonts, amsmath}
\usepackage{epsfig}
\usepackage{graphics, latexsym,amssymb}
 \def\caixavazia #1#2{{\dimen0=#1 \advance\dimen0 by -#2       
       \dimen1=#1 \advance\dimen1 by #2                       
        \vrule height #1 depth #2 width #2                    
        \vrule height 0pt depth #2 width #1                   
        \llap{\vrule height #1 depth -\dimen0 width \dimen1}% 
       \hskip -#2                                             
       \vrule height #1 depth #2 width #2}}                   
 \def\caixaqed{\mathord{\caixavazia{6pt}{.4pt}}\;}

\def\endpf{\hfill $\caixaqed$}
\def\qed{\hfill $\caixaqed$}
\def\S3{{\mathbb S^3}}
\def\R3{{\mathbb R^3}}
\def\proof{\noindent \bf Proof. \rm}

\font\teneufm=eufm10
\font\seveneufm=eufm7
\font\fiveeufm=eufm5
\newfam\eufmfam
\textfont\eufmfam=\teneufm
\scriptfont\eufmfam=\seveneufm
\scriptscriptfont\eufmfam=\fiveeufm

\newfam\msbfam
\font\tenmsb=msbm10  \textfont\msbfam=\tenmsb
\font\sevenmsb=msbm7  \scriptfont\msbfam=\sevenmsb
\font\fivemsb=msbm5    \scriptscriptfont\msbfam=\fivemsb

\newtheorem{theorem}{Theorem}[section]

\newtheorem{lemma}{Lemma}[section]
\newtheorem{definition}{Definition}[section]
\newtheorem{remark}{Remark}[section]
\newtheorem{corollary}{Corollary}[section]

\begin{document}

\title{\bf Convergent sequences of closed minimal surfaces embedded in $\S3$}

\author{Fernando A. A. Pimentel}

\maketitle

\begin{quotation}\noindent
{\bf Abstract:} \sl given two minimal surfaces embedded in $\S3$ of genus $g$
we prove  the existence of a sequence of non-congruent compact minimal surfaces embedded in $\S3$ of genus $g$ that converges in $C^{2,\alpha}$ to a compact embedded minimal surface provided some conditions are satisfied.  These conditions also imply that, if any of these two surfaces is embedded by the first eigenvalue, so is the other.
\end{quotation}

\allowdisplaybreaks

\section{Introduction}

The objective of this paper is to establish sufficient conditions for the existence of  sequences of non-congruent compact minimal surfaces embedded in $\S3$. These sequences shall have  special features; among them,   convergence in $C^{2,\alpha}$ to a compact embedded minimal surface. 

Similar conditions for the existence of sequences of minimal surfaces were given in the fourth section of \cite{p1} as a necessary step to prove the well known Lawson conjecture for minimal tori 
embedded in $\S3$. As in this previous work, our results here are a means to an end, since they are intended to be used in an investigation in $\S3$ of the Yau's conjecture for the first eigenvalue of the laplacian of minimal embeddings.

Although we, in a sense,  generalize and adapt results and ideas in \cite{p1} to the peculiarities of the present investigation, it is evident  that a mere generalization of the results in \cite{p1} to compact embedded minimal surfaces of arbitrary genus, though quite straightforward, does not provide what is needed to attain our goals.  In fact, the sequences that we deal with here are different from the sequences that are studied in \cite{p1}, which leads to important differences from this work to section 4 of \cite{p1}. Thus   a notion of distance from a minimal torus to the Clifford torus is explicitly controlled in \cite{p1}. In this paper,  we   discuss a different approach to define sequences of compact embedded minimal surfaces, that is, by controlling their first eigenvalue. Each minimal surface in a sequence is then characterized by being embedded by a different first eigenvalue (assumed the previous existence of two such surfaces embedded by distinct first eigenvalues.)

As said above, the existence of these sequences is not without consequences. Indeed, if we assume that the conditions for their existence (that are very loose) are satisfied by any pair of non-congruent compact minimal surfaces of a given genus $g$ embedded in $\S3$, we can affirm that either such surfaces are never isolated in $C^{2,\alpha}$ or that there exists, up to congruences, only one compact surface of genus $g$ embedded in $\S3$.

Another consequence of  the assumption that the conditions of Theorem \ref{tp2} are satisfied by any minimal surface embedded in $\S3$ is a proof of the well known Yau's conjecture for minimal surfaces in $\S3$ (see \cite[prob. 100]{y}.) A key ingredient for proving this result, besides Theorem \ref{tp2} and the fulfillment of its conditions by arbitrary minimal surfaces, is the fact that all Lawson's embedded minimal surfaces (the embedded examples in \cite{l1}) are embedded by the first eingenvalue, proved recently by Choe-Soret (see \cite{cs}.) In our Final Remarks, we discuss how to complete a proof of Yau's conjecture.

We complement  these brief comments on the nature of our results by fully stating them below (Section \ref{sst}). Next, we discuss the  organization of this paper (Section \ref{sorg}). Some questions concerning our previous work referenced  in \cite{p1} are answered here. We pose these questions and point  to where they are discussed in Section \ref{sorg}.

\section{Statement of the main results} \label{sst}

We use here the terminology and notation established in \cite[Section 2]{p1}. Thus, we let $\langle\hspace{.02in},\hspace{.02in}\rangle$ be the inner product of the euclidean space $\mathbb R^4$ and also let $||\ ||$ denote its associated norm so $\S3=\{p\in\mathbb R^4\hspace{.02in}:\hspace{.02in}||\hspace{.01in}p\hspace{.02in}||=1\}$. Given a unit vector $v\in \mathbb R^4$ we have the equator (big sphere) $S(v)=\{p\in \S3 \hspace{.02in}|\hspace{.02in}\langle v,p\rangle=0\}$

For what follows, we refer to \cite{sp} for the definitions and properties of the curvature of a curve  and the mean curvature of a surface in $\S3$.

Along this work, we assume the existence of non-congruent compact minimal surfaces of genus $g$, namely $N_1$ and $N_2$, embedded in $\S3$. We also assume that there exists $\lambda_1:S^1 \to N_1$ and $\lambda_2:S^1 \to N_2$, embedded closed real analytic curves in the intersection of $N_1$ and $N_2$, respectively, with equators of $\S3$, that are parameterized by arc length. These parameterizations may be chosen in such a way that $\lambda_1(0)$ and $\lambda_2(0)$ are points where the curvature attains its maximal value at the respective curve.  

Recall the variational characterization of the first (positive) eigenvalue  of a compact surface $M$:

\begin{equation}\label{epa}
\kappa_1(M)=\inf\left\{\dfrac{\int_{M}|\hbox{grad}_{M} f|^2\ dM}{\int_M f^2 dM}\hspace{.06in}\Big|\hspace{.06in}{f\in C^1(M)}\hbox{ and }{\int_M f dM=0}\right\}.
\end{equation}

It is well known that the coordinate functions provided by the immersion of $\S3$ into $\mathbb R^4$, when restricted to a minimal surface in $\S3$, are eigenfunctions of the laplacian of the surface, with eigenvalue equal to two (Takahashi \cite{tk}). Besides, two compact minimal surfaces embedded in $\S3$ are diffeomorphic by a diffeomorphism of $\S3$ (\cite[Th. 3]{l2}). Thus it seems reasonable to impose the conditions below on $N_1$ and $N_2$:

\begin{enumerate}
\item the surface  $N_2$ and only the surface $N_2$ is embedded by the first eigenvalue;
\item the curvatures of $\lambda_1$ at $\lambda_1(0)$ and of $\lambda_2$ at $\lambda_2(0)$ are equal;
\item there exists a $C^{\infty}$ diffeomorphism $X_{12}$ of $\S3$ taking $N_1$ to $N_2$, $\lambda_1$ to $\lambda_2$ and $\lambda_1(0)$ to $\lambda_2(0)$ whose restriction to $\lambda_1$ preserves orientation for any previously chosen orientations of $\lambda_1$ and $\lambda_2$; 
\end{enumerate}
 (we recall that we orient a regular curve by choosing  a unit tangent vector field along it). Thus there holds that
 
 \begin{theorem}\label{tp1} Assume that there exists $N_1$ and $N_2$,   noncongruent compact minimal surfaces of genus $g$ embedded in $\S3$ as above, that also satisfy conditions (2) and (3). Then, fixed $0<\alpha<1$, there exists a sequence of noncongruent compact minimal surfaces of genus $g$ embedded in $\S3$ converging to $N_2$ in $C^{2,\alpha}$.
 \end{theorem}

However, if we do not drop condition (1) from our hypotheses on $N_1$ and $N_2$, we also have the ensuing theorem. 
 
 \begin{theorem}\label{tp2} If there exists minimal surfaces $N_1$ and $N_2$ as above satisfying conditions (1), (2), and (3) then, fixed $0<\alpha<1$,  there exists a sequence of compact embedded minimal surfaces $(M_n)$ of genus $g$ such that 
\begin{itemize} 
\item[(i)] the sequence $(M_n)$ converges  in $C^{2,\alpha}$ to a compact minimal surface $M$ of genus $g$ (not necessarily congruent to $N_2$) that is embedded in $\S3$ by the first eigenvalue;
\item[(ii)] if $n>m$ then the first eigenvalue of the laplacian on $M_n$ is greater than the first eigenvalue of the laplacian on $M_m$.
\end{itemize}
 \end{theorem}

 \section{Some remarks on the organization  and contents of this paper}\label{sorg}
 
 \noindent
{\bf A) Structure of the paper}\bigskip
 
We  study in Section \ref{s1} the properties of the functions $\tau^\alpha$ and $\tau_*^\alpha$, $0<\alpha\leq 1$. These functions have the same roles since they  measure  a distance from a $C^{2,1}$ diffeomorphism of $\S3$ to the identity map. The functional $\tau^\alpha$  was defined in \cite{p1}  and  is used here in the proof of Theorem \ref{tp1}, which generalizes the similar result in \cite{p1}. From the functional $\tau^\alpha$ we define and characterize the new functional  $\tau_*^\alpha$ on the set of $C^{2,\alpha}$ diffeomorphisms of $\S3$ in order to prove Theorem \ref{tp2}.

The functions $\xi_{n,\alpha}^{t,k}$ and $\Phi_{n,\alpha}$\hspace{.01in} (see definitions \ref{dxi} and \ref{dphi}, respectively), versions of their namesakes in \cite{p1}  for use in the proof of Theorem \ref{tp2}, are presented in Section \ref{s2}. We  provide below a brief description of these functions.

The function $\xi_{n,\alpha}^{t,k}$ (Definition \ref{dxi}) assigns to a $C^{2,1}$ closed  curve $\lambda$ embedded in an equator of $\S3$ the infimum in a set   of closed surfaces of genus $g$ embedded in $\S3$ containing the curve $\lambda$ a function of  the first eigenvalue  of such surfaces. These surfaces, among their features, also have mean curvature equal to zero at prescribed  points and are taken to $N_2$, a fixed compact surface minimally embedded in $\S3$, by a diffeomorphism $X$ of $\S3$ such that $\tau_*^1(X)$ is majored by the constant $k>0$. 

The value of $\Phi_{n,\alpha}(t)$ (Definition \ref{dxi}) is then defined in Section \ref{s2} as the infimum of $\xi_{n,\alpha}^{t,\kappa}$ in a set of embedded closed curves in equators of $\S3$ containing an arc of curve denoted $\lambda_1^t$. The arc $\lambda_1^t$, in turn, is contained in a closed curve in $N_1$, the  compact minimal surface in $\S3$ diffeomorphic (though not congruent) to $N_2$, where $N_1$ and $N_2$ are minimal surfaces defined in Section \ref{sst} whose existence we will assume.

We will verify that the function $\Phi_{n,\alpha}$ is continuous in Section \ref{scont}. This result will enable us to produce in Section \ref{spmr} a sequence of embedded minimal surfaces satisfying the conditions in the statement of theorems \ref{tp1} and \ref{tp2}, proving thus these theorems.

In our Final remarks, we consider how to use theorems \ref{tp1} and \ref{tp2} to investigate some open questions concerning minimal surfaces of arbitrary genus embedded in $\S3$.

Finally, for the sake of the fluidity of the exposition, we collect in appendices \ref{A1} and \ref{A2} two technical results, lemmas \ref{laux} and \ref{lfal}, that are used in the proof of lemmas \ref{lcon1} and 
\ref{lcon2}, respectively.\bigskip

\noindent
{\bf B) Sketch of the proof of the main theorems}\bigskip

The process presented in \cite{p1} and further developed here, used to prove our main results,  may be summarized as follows (we particularize our discussion to the proof of Theorem \ref{tp2}, because the proof of Theorem \ref{tp1} is similar.) Assume the conditions of Theorem \ref{tp2} satisfied (thus, unlike $N_2$, the surface $N_1$ is not embedded by the first eigenvalue.) We then proceed by the ensuing steps: \bigskip

\noindent
\hspace{.05in}STEP 1: we produce sets of  $n$ points $V_n\subset N_2$ such that $V_n\subset V_m$ if $n\leq m$ and $\bigcup_{n=1}^\infty V_n$ is dense  in $N_2$.\medskip

\noindent
\hspace{.05in}STEP 2: we prove: (i) that the functions  $\Phi_{n,\alpha}$ are uniformly bounded by a constant $C$ that does not depend on $n$ nor on $\alpha$; (ii) that these functions are continuous.\medskip

\noindent
\hspace{.05in}STEP 3: let $\gamma\in [\kappa_1,2]$, where $\kappa_1$ is the first eingenvalue of $N_1$. We obtain that the continuity of the functions $\Phi_{n,\alpha}$ implies the existence of a diffeomorphism  of $\S3$ denoted by $X_{n,\gamma}$ such that the first eigenvalue of $X_{n,\gamma}(N_2)$ is equal to $\gamma$ and, if $p\in V_n$, then the mean curvature of $X_{n,\gamma}(N_2)$ is equal to zero at $X_{n,\gamma}(p)$.\medskip

\noindent
\hspace{.05in}STEP 4: we prove that a subsequence $X_{n_j,\gamma}$ converges in $C^{2,\alpha}$ to a diffeomorphism of $\S3$, say $X_\gamma$. We will also have the convergence of the inverse functions $X_{n_j,\gamma}^{-1}$ to $X_{\gamma}^{-1}$ (the existence of these convergent subsequences is a consequence of Lemma \ref{lbas} and the uniform boundness  of  the values of $\tau_*^\alpha$ along these sequences. These results are obtained in sections \ref{s1} and \ref{s2}, respectively.) 

\noindent
OBS. 1: the set $X_\gamma(\bigcup_{n=1}^\infty V_n)$ shall be dense in $X_\gamma(N_2)$,
because the restriction of $X_\gamma$ to $N_2$ is a diffeomorphism into $X_\gamma(N_2)$ and diffeomorphisms take dense sets to dense sets.

\noindent
OBS. 2: since points in $X_{n_j,\gamma}(V_k)$ have mean curvature equal to zero in the surface $X_{n_j,\gamma}(N_2)$ whenever
$n_j>k$, the $C^{2,\alpha}$ convergence $X_{n_j,\gamma}\to X_\gamma$, which  implies that $X_{n_j,\gamma}(V_k)$ converges to
$X_\gamma(V_k)$, also implies that points in $X_\gamma(V_k)$ have mean curvature zero in $X_\gamma(N_2)$. This fact along with OBS. 1 implies that $X_\gamma(N_2)$ has mean curvature equal to zero at every point.

\noindent
OBS. 3: The $C^{2,\alpha}$ convergence $X_{n_j,\gamma}\to X_\gamma$ also implies that the first eigenvalue of $X_\gamma(N_2)$ is equal to $\gamma$. \medskip

\noindent
\hspace{.05in}STEP 5: we prove by Lemma \ref{lbas} that the diffeomorphisms $X_\gamma$ are $C^{2,1}$ diffeomorphisms at which the function $\tau^1_*$ is bounded. \medskip

\noindent
\hspace{.05in}STEP 6: proceeding as in STEP 4, by STEP 5 and Lemma \ref{lbas} there exists  an increasing sequence  $\gamma_j \to 2$ such that the sequence $(X_{\gamma_j})$ converges in $C^{2,\alpha}$ to a diffeomorphism $X$ of $\S3$ such that $X(N_2)$ is minimal and $X(N_2)$ is embedded by its first eigenvalue (which is two by a well known Theorem by Takahashi in \cite{tk}).\bigskip

\noindent
{\bf C) Answering some questions regarding our previous work \cite{p1}}\bigskip 

Along this paper, we will take  advantage of the opportunity to answer
some questions on our previous work \cite{p1} on the steps above.

\begin{itemize}
\item A recurrent argument (used, e.g., in steps 4 and 6 above) goes as follows: a sequence $(X_n)$ of $C^{2,1}$ diffeomorphisms of $\S3$  converges in $C^{2,\alpha}$, $0<\alpha<1$, to $X\in C^{2,\alpha}$, a diffeomorphism of $\S3$ itself.  It is also assumed that $\tau_*^1(X_n)\leq k$ $\forall\hspace{.02in} n$ for some positive constant $k$. It is clear that $\tau_*^1$ is not continuous in the $C^{2,\alpha}$ topology (in fact, it is not defined in this topology.) Anyway, some extra information is recovered from these conditions on the sequence $(X_n)$: it is contained in $C^{2,1}\cap C^{2,\alpha}$ and $\tau_*^1$ is bounded in it. From these special conditions we prove in Lemma \ref{lbas} that $X$ is also in $C^{2,1}$ and that $\tau_*^1(X)\leq k$. This is thus a case where regularity is recovered in special conditions. We ask the reader to study Lemma \ref{lbas} and its proof, then confer Remark \ref{resc} to settle this issue. We also ask the reader to look at Remark \ref{rred}.
\item In Step 4, a sequence $(M_{n_j})$ converges in $C^{2,\alpha}$ to a surface $M$. The surfaces $M_{n_j}$ have an increasing number of points with mean curvature zero. The argument requires that these points converge to a dense set of points in $M$. This issue is treated in observations 1 and 2 to Step 4 above and fully discussed in the proof of Theorem \ref{tp2} (Section \ref{spmr}.)
\item A critical point in Step 2 is the proof that $\lim_{t\to 1^-}\Phi_{n,\alpha}=0$. This limit is explicitly proven in Lemma \ref{lcon2} with a necessary clarification of an important point provided by Lemma \ref{lfal}.
\end{itemize}

\vspace{.15in}

\section{The functions $\tau^\alpha$ and $\tau^\alpha_*$: measuring distances between diffeomorphisms}\label{s1}

\vspace{.1in}

\noindent
{\bf A) The Function $\tau^\alpha$}

\vspace{.2in}

The function $\tau^\alpha$ was introduced in \cite{p1} and shall be used in the proof of Theorem \ref{tp1}. We will define it again and restate its main properties. 
As in \cite[Sect. 4]{p1} we let $A_2=\{x\in \mathbb R^4:\hspace{.03in} 1/2<|x|<2\}$ and have then for $C^{2,\alpha}$ maps $F:A_2\to A_2$  the usual Holder norms and seminorms for a subset of a euclidean space. Let $f:A_2\to\mathbb R$ be a $C^{0,\alpha}$ function. If $F(x)=(F_1(x),F_2(x),F_3(x),F_4(x))$ we let
\begin{eqnarray}\label{eqnh1}
[f]_{\alpha}&=&\sup\left\{\frac{|f(x)-f(y)|}{||x-y||^\alpha}\hspace{.06in}\big|\hspace{.06in}{x,y\in A_2},\ {x\neq y}\right\},\\ \nonumber
||F||_{C^{2,\alpha}}&=& \sum_{i=1}^4\sup_{A_2}|F_i(x)|+ \sum_{1\leq i,j \leq 4} \sup_{A_2}\left|{D^j F_i(x)}\right|\\ \label{eqnh2}
&&+\sum_{1\leq j\leq k\leq 4}\left[\left(\sum_{i=1}^4\sup_{A_2}\left|D^{jk} F_i(x)\right|
\right)+\left( \sum_{i=1}^4[D^{jk} F_i]_{\alpha}\right)\right],
\end{eqnarray}
where $||\hspace{.02in}\hspace{.02in}||$ is the euclidean norm in $\mathbb R^4$.
\medskip

\noindent
{\bf Terminology:} a $C^{2,\alpha}$ diffeomorphism $\xi$ of $\S3$ is canonically 
extended to a $C^{2,\alpha}$ diffeomorphism $X$ of $A_2$: for all $v\in \S3$ and $1/2< r<2$  
we let $X(rv)=r\xi(v)$. We will then say that $X$ is a diffeomorphism of $\S3$ canonically extended to $A_2$.\medskip

\begin{definition}\label{dxcal} Let $\mathcal{X}^\alpha$ be the set of  $C^{2,\alpha}$ diffeomorphisms of  $\S3$ canonically extended to $A_2$.
\end{definition}

The introduction of the spaces $\mathcal X^\alpha$ in \cite{p1} was meant to  make use of the vector space 
structure of $\mathbb R^4$ along with the  holder norms and the canonical embedding theorems for $C^{k,\alpha}$ functions on euclidean spaces in the most straightforward manner. Since the diffeomorphisms of $\mathcal X^\alpha$ are diffeomorphisms of $\S3$ canonically extended to $A_2$ in the sense above, working in $\mathcal X^\alpha$ is equivalent to work with $C^{2,\alpha}$ diffeomorphisms of $\S3$.

For the sake of completeness we restate here \cite[Lem 4.1]{p1}, an immediate consequence of \cite[Rem. 4.1]{p1} and \cite[Th. 3.1]{ad}:

\begin{lemma}\label{lemb}  If  $0<\nu<\lambda\leq 1$  then the embedding $\mathcal{X}^\lambda\to \mathcal{X}^\nu$ exists and is compact. 
\end{lemma}

Following \cite{p1},  we let  $I:A_2\to A_2$ be the identity map and define the following functional on the diffeomorphisms of $A_2$:
\begin{equation}
\tau^\alpha(X)= ||X-I||_{C^{2,\alpha}}+||X^{-1}-I||_{C^{2,\alpha}},
\end{equation}
where $0<\alpha\leq 1$.

The following lemma is identical to \cite[Lem. 4.2]{p1}:

\begin{lemma}\label{lbas}  If $0<\alpha\leq 1$ let $(X_k)$ be a sequence of $C^{2,\alpha}$ diffeomorphisms of $\S3$ canonically extended to $A_2$ as above. Assume that $\tau^\alpha(X_k)\leq C<\infty\ \forall\hspace{.02in}k\in \mathbb N$.  Let  $0<\beta<\alpha$. Then 
\begin{itemize}
\item[(i)] a subsequence of $(X_k)$ converges  in
$C^{2,\beta}$ to a diffeomorphism $X$ of $\S3$ canonically extended to $A_2$; 
\item[(ii)] besides being a map in $C^{2,\beta}$, the map $X$ is also a $C^{2,\alpha}$ diffeomorphism and ${\tau^{\alpha}(X)}\leq C$. 
\item[(iii)] if $X_k \to X$ in $C^{2,\alpha}$, then $X$ is a $C^{2,\alpha}$ diffeomorphism of $\S3$ canonically extended to $A_2$, and there exists $\lim_{k\to\infty}\tau^\alpha(X_k)=\tau^\alpha(X)$.
\end{itemize}
\end{lemma}

\noindent
{\proof} 
See the proof of \cite[Lemma 4.2]{p1}.

\begin{remark}\label{resc}
Lemma \ref{lbas} is a direct consequence of the compact embedding results  for Holder spaces (\cite[Th. 3.1]{ad}) and is meant to provide compactness criteria  for sequences of diffeomorphisms of $A_2$. Thus what it is affirmed  in Lemma \ref{lbas} is that the space of $C^{2,\alpha}$ diffeomorphisms of $A_2$ is such that any sequence in this space is precompact in the space of $C^{2,\beta}$ diffeomorphisms of $A_2$, provided that $0<\beta<\alpha\leq 1$ and that $\tau^\alpha$ is bounded in the sequence. 

We must also point out that Lemma \ref{lbas}  is a typical instance of a result that characterizes situations where, under certain conditions, regularity is recovered. Thus a sequence in $C^{2,\alpha}$ with $\tau^\alpha$ bounded by $C$ has a subsequence that converges in $C^{2,\beta}$ to a $C^{2,\beta}$ diffeomorphism  $X$ of $A_2$ that also is in $C^{2,\alpha}$. Moreover, Lemma \ref{lbas} also guarantees that  $\tau^\alpha(X)\leq C$ (although we don't know its exact value).  

Therefore, if $S^\circ$ is the closure in $\mathcal{X}^\beta$ of a subset $S$ of $\mathcal{X}^\alpha$ on which $\tau^\alpha$ is bounded by $C$  then $S^\circ$ is contained in $\mathcal{X}^\alpha$. Besides, the restriction of $\tau^\alpha$ to $S^\circ$ is bounded by $C$ as well. 

We apologize for these redundant remarks but we  think that they are necessary to avoid misunderstandings.
\end{remark}

We now define a distance in $\mathcal X^\alpha$ from  which the function $\tau^\alpha$ is a  particularization:

\begin{definition}\label{ddist} Given $X,$ $Y\in X^\alpha$ let 
$$d^{\alpha}(X,Y)=||X-Y||_{C^{2,\alpha}}+||X^{-1}-Y^{-1}||_{C^{2,\alpha}}.$$
\end{definition}
Hence, $\tau^\alpha(X)=d^\alpha(X,I)$ and a  result similar to Lemma \ref{lbas} holds for the binary application $d^\alpha :\mathcal X^\alpha\times \mathcal X^\alpha \to \mathbb R$. 

\vspace{.3in}

\noindent
{\bf B) The function $\tau^\alpha_*$}

\vspace{.2in}

We shall identify diffeomorphisms  which differ by a congruence of $\S3$: given $X$, $Y\in \mathcal \mathcal X^\alpha$, then $X\sim Y$ if $X\circ Y^{-1}\in O(4)$. Let $\bar X$ be the equivalence class of $X$ in $\mathcal X^\alpha/\sim$ and $\pi:\mathcal X^\alpha\to \mathcal X^\alpha/\sim$ the canonical projection (hence, $\pi(X)=\bar X$). 
The distance $d^\alpha_*$ in $\mathcal X^\alpha/\sim$ is given by 
\begin{equation}\label{edd}
d^\alpha_*(\bar X,\bar Y)=\inf\{\hspace{.02in}d^\alpha(X,Y)\ :\  \pi(X)\in \bar X,\ \pi(Y)\in \bar Y\},
\end{equation}

We define the function $\tau^\alpha_*$ accordingly: 
\begin{definition}
 If $0<\alpha\leq 1$ and $X\in\mathcal X^\alpha$ let
$\tau^\alpha_*(X)= d_*^\alpha(\bar X,\bar I)$. 
\end{definition}

\begin{remark}\label{rbas} We readily see that the conclusions of Lemma \ref{lbas} and the reasoning in  Remark \ref{resc} also apply to the function $\tau^\alpha_*$.
\end{remark}

The functions $\tau^\alpha_*$ then identify congruent diffeomorphisms of $\S3$ and have the good  features of the $\tau^\alpha$ functions, so their usefulness. This function shall be used in the  definition of the functions $\xi_{n,\alpha}^{t,k}$ for the proof of Theorem \ref{tp2}.

\section{ The functions $\xi_{n,\alpha}^{t,k}$ and $\Phi_{n,\alpha}$}\label{s2}

As in \cite{p1},
we will produce a set that is dense  in a given $n$-torus. Such dense set shall be obtained from the union of subsets  with an increasing number of more or less equally distributed points. In \cite{p1}, a dense subset of the Clifford  torus was defined  by  reticulating the Clifford torus by its longitudes and latitudes. An arbitrary $n$-torus does not have longitudes and latitudes and cannot be reticulated in this way. But it can always  be triangulated. 

We recall that a two dimensional triangulation  $T$ is a finite subset of $\mathbb R^2$ whose elements are triangles (i.e., the interior of triangles), open segments and points such that any two elements of $T$ are disjoint and $T$ contains both vertices of any of its segments and all the sides and vertices of any of its triangles. A polyhedron is  the body of a triangulation (see \cite{al} for definitions and results on triangulations.)

From the previous section we recall the compact minimal surface $N_2$. It is well known that $N_2$ is homeomorphic to a finite polyhedron $P$ (see \cite[Remark, p. 74]{al}), i.e., a polyhedron whose triangulation has a finite number of triangles, segments and points. We also observe that all points in the triangulation of $P$ are vertices of some of its triangles. The triangulation of $P$ then induces  
a triangulation $T_0$ of $N_2$ via the homeomorphism: we call vertices, segments and triangles the image in $N_2$ by the homeomorphism of the vertices, segments and triangles of the polyhedron $P$.

 Successive  barycentric subdivisions (see \cite[p. 78]{al})  of the triangulation of the polyhedron $P$ then induce successive subdivisions of $T_0$, which, following \cite{al}, we also call barycentric subdivisions. Let then $T_n$ be the $n^{th}$ barycentric subdivision  of $T_0$.

\begin{definition} \label{dvn} Let $V_n$ be the set of vertices of $T_n$.
\end{definition}

We recall that dense sets are taken to dense sets by homeomorphisms.  Hence, $\bigcup_{n=0}^{\infty} V^n$ is a countable dense subset of $N_2$ as desired. Moreover, $V^m\subset V^n$ if $m\leq n$.

We hereon define some spaces of diffeomorphisms of $A_2$ as well as functions on these spaces.

\begin{definition}\label{d0} Let $\lambda_2\subset N_2$ be the closed curve in $N_2$ defined in section \ref{sst}. 
Let $\lambda$ be a closed $C^{2,1}$ curve embedded in a equator. Let then $\Omega_{n}^\lambda$ be the set of  $C^{2,1}$ diffeomorphisms $X$ of $\S3$ canonically extended to $A_2$ such that both $X(\lambda_2)$ is congruent to $\lambda$  as pointsets and the the compact surface  $X(N_2)$
has mean curvature equal to zero at the points in $X(V_n)$.
\end{definition}

We stress that the mean curvature of $X(V_n)$ (and of any surface considered here) is relative to its  immersion in $\S3$.

Proceeding  as in \cite[Remark 4.2]{p1}, we verify that the sets $\Omega_{n}^\lambda$ are not empty. In fact, given a diffeomorphism $X$ of $\S3$, we can prescribe the mean curvature of  $X(N_2)$ at the image of a discrete set of points (like $V_n$) by perturbing the diffeomorphism  $X$ in neighborhoods of these points.

\begin{definition} \label{domega2}
If $0<\alpha<1$ we let
$\Omega_{n,\alpha}^{\lambda,{k}}$ be the closure in $C^{2,\alpha}$ of the set of diffeomorphisms $X\in \Omega_{n}^\lambda$ such that $\tau_*^1(X)<{k}$.
\end{definition}

\begin{remark}\label{rred} We observe that $\Omega_{n}^\lambda$ is in $C^{1,1}$ while $\Omega_{n,\alpha}^{\lambda, k}$ is in  $C^{1,\alpha}$.  However, from Lemma \ref{lbas} (see also remarks \ref{resc} and \ref{rbas}) we also have that $\Omega_{n,\alpha}^{\lambda,{k}}$ is a set of $C^{2,1}$ diffeomorphisms of $\S3$ and thus 
$\Omega_{n,\alpha}^{\lambda,{k}}\subset \Omega_{n}^\lambda$ as well.
\end{remark}

\begin{definition}\label{dsigma} For $X\in \mathcal X^\alpha$, $0<\alpha\leq 1$, we let $\sigma(X)$ be the first nonzero eigenvalue of $X(N_2)$, seen as a compact surface embedded in $\S3$. 
\end{definition}
Thus it is readily seen that $\sigma$ is continuous in $\mathcal X^\alpha$ (see Lemma \ref{lsigma}).

\smallskip
\noindent{\bf Notation:} {\it we denote by $\kappa_1$, $\kappa_2$ the first nonzero eigenvalues of the embeddings of $N_1$, $N_2$ into $\S3$, respectively.}\smallskip

We now define a function on the set of simple $C^{2,1}$ curves in an equator of $\S3$.

\begin{definition}\label{dxi}
 If $\Omega^{\lambda,k}_{n,\alpha}$ is not empty we set for $0\leq t\leq 1$:
\begin{equation}\label{eqxi}
\xi_{n,\alpha}^{t,k}(\lambda)=\inf\{\hspace{.02in}|\sigma(X)-(t\kappa_2+(1-t)\kappa_1)|\hspace{.03in}:\hspace{.03in}{X\in \Omega^{\lambda,k}_{n,\alpha}} \}
\end{equation}
\end{definition} \medskip

We below set the notation for definitions \ref{dl} and \ref{dphi}.

\begin{itemize}\label{prod}
\item Let $\mu,$ $\nu: [0,1]\to\S3$ be $C^{2,1}$ curves in a same equator $S(v)$ of $\S3$ with $\mu(0)=\nu(1)$ and $\mu(1)=\nu(0)$. If $\mu\cup\nu$ is the trace of a $C^{2,1}$ curve in $S(v)$ we denote by the product notation $\mu\nu$ any arclength parameterization of this curve.
\item If $\lambda$ is an  arclength parameterized closed curve in $S(v)$ with length $L$,  we denote by $\lambda^t$ the restriction of $\lambda$ to $[-(1-t)L/2 ,(1-t)L/2]$.
\item We assume that the conditions of Theorem \ref{tp1} are satisfied and, among them, we recall the existence and properties of the diffeomorphism of $\S3$ denoted by $X_{12}$, defined in Section \ref{sst}. We consider then the canonical extension of $X_{12}$ to $A_2$, and call it again $X_{12}$, abusing of the notation. Thus $X_{12}^{-1}\in \Omega_{n}^{\lambda_1}$ for all $n\in \mathbb N$; Obviously we have that $\tau_*^\alpha(X_{12})=\tau_*^\alpha(X_{12}^{-1})$. 
\end{itemize}

\begin{definition}\label{dl}  Let $\kappa >\tau_*^\alpha(X_{12})$. If $0<\alpha,\ t<1$ let
 $\mathcal{S}_{n,\alpha}^t$ be the set of $C^{2,1}$ arcs $\mu:[0,1]\to S(v)$ in the equator $S(v)$ containing $\lambda_1$, with the same extremities as $\lambda^t_1$, such that 
\begin{itemize}
\item[(1)] the closed curve $\lambda_1^t\mu$ is a $C^{2,1}$ curve embedded in $S(v)$; 
\item[(2)] the set $\Omega^{\lambda^t_1\mu,\kappa}_{n,\alpha}$ is defined and not empty.
\end{itemize}
\end{definition}

 The set $\mathcal{S}_{n,\alpha}^t$ was also defined in \cite{p1}, where it was observed that $\mathcal{S}_{n,\alpha}^t$ is not empty because it contains at least the curve $\mu$ whose trace is the closure of
$\lambda_1\backslash \lambda^t_1$. Hence, $\lambda_1^t\mu$ is a reparameterization of $\lambda_1$  (observe that
that $\tau_*^\alpha(X_{12})\geq \inf \{\tau_*^\alpha(X)\hspace{.02in}:\hspace{.02in} X\in\Omega_{n}^{\lambda_1}\}$).

\begin{definition}\label{dphi} For $t\in (0,1)$ let
\begin{equation}\label{eqphi}
\Phi_{n,\alpha}(t)=\inf_{\mu\in \mathcal{S}_{n,\alpha}^t}\xi_{n,\alpha}^{t}(\lambda_1^t\mu)
\end{equation}
and let $\Phi_{n,\alpha}(0)=\Phi_{n,\alpha}(1)=0$.
\end{definition}

\section{Some continuity results}\label{scont}

We now discuss the continuity of the functions $\sigma$ and $\Phi_{n,\alpha}$ (see definitions \ref{dsigma} and \ref{dphi}). The proof of Lemma \ref{lsigma} below is straightforward.

\begin{lemma}\label{lsigma} If $0<\alpha\leq 1$ then the function $\sigma$  is continuous in $\mathcal X^\alpha$.
\end{lemma}

\begin{lemma}\label{lcon0} Fixed $0<\alpha<1$, for every $n\in \mathbb N$ each function $\Phi_{n,\alpha}$ is bounded by a same constant $C$ that does not depend on $n$ nor on $\alpha$.  
\end{lemma} 

\proof 
We observe that, since $N_1=X_{12}^{-1}(N_2)$ is minimal, we have that $X_{12}^{-1}\in \Omega_{n,\alpha}^{\lambda_1,\kappa}$ and thus the closure of $\lambda_1\backslash\lambda_1^t$ belongs to $\mathcal{S}_{n,\alpha}^t$ for every $n\in \mathbb N$. Hence, $$\Phi_{n,\alpha}(t)\leq \xi^{t,\kappa}_{n,\alpha}(\lambda_1)\leq|\sigma(X_{12}^{-1})-(t\hspace{.02in}\kappa_2+(1-t)\kappa_1)|
=(1-t)(\kappa_1+\kappa_2).$$
Let thus $C=\kappa_1+\kappa_2$. \endpf

\begin{remark}\label{rconvc} We will discuss below the convergence of sequences of embedded closed curves denoted by the product notation described above. In order to normalize the parameterizations, we will assume that any curve  $\lambda_1^t\mu$ obtained by the product of $\lambda_1^t$ with a suitable arc is an application from $S^1$ to an equator of $\S3$ parameterized by arc length that is oriented as $\lambda_1^t$ with $\lambda_1^t\mu(0)=\lambda_1(0)$.
\end{remark}
\medskip

\begin{lemma}\label{lcon1} For every $n\in \mathbb N$ and $\alpha\in (0,1)$ the function $\Phi_{n,\alpha}$ is continuous in the open interval $(0,1)$.
\end{lemma}

\proof Assume that either $ \nexists\lim_{t\to t_0}\Phi_{n,\alpha}(t)$ or $\lim_{t\to t_0}\Phi_{n,\alpha}(t)\neq\Phi_{n,\alpha}(t_0)$
for some $0<t_0<1$. Since $\Phi_{n,\alpha}$ is bounded, these assumptions are equivalent to the existence of a sequence $\theta_k\to t_0$ such that  $\exists \lim_{k\to\infty}\Phi_{n,\alpha}(\theta_k)\neq\Phi_{n,\alpha}(t_0).$ 
We then divide our analysis in two cases.\medskip

\noindent
{\bf Case 1:} assume that $\lim_{k\to\infty}\Phi_{n,\alpha}(\theta_k)>\Phi_{n,\alpha}(t_0).$
Let then 

\begin{equation}\label{eep}
2\epsilon=\lim_{k\to \infty}\Phi_{n,\alpha}(\theta_k)-\Phi_{n,\alpha}(t_0).
\end{equation}

 We here follow closely the proof of \cite[Lem 4.3]{p1}. Thus, 
we let $\mu_j^{t_0}\subset \mathcal{S}_{n,\alpha}^{t_0}$ be a sequence  of not necessarily distinct arcs $\mu_j^{t_0}$ such that $\lim_{j\to\infty}\xi_{n,\alpha}^{t_0}(\lambda_1^{t_0}\mu_j^{t_0})=\Phi_{n,\alpha}(t_0)$. 

If $j$ is sufficiently large then we may assume that $$0\leq \xi_{n,\alpha}^{t_0} (\lambda_1^{t_0}\mu_j^{t_0})-\Phi_{n,\alpha}(t_0)<\epsilon.$$ Thus there exists 
$X\in \Omega^{\lambda_1^{t_0}\mu^{t_0}_j,\kappa}_{n,\alpha}$ such that $\tau_*^1(X)<\kappa$ 
(by the definition of the sets $\Omega^{\lambda_1^{t_0}\mu^{t_0}_j,\kappa}_{n,\alpha}$ and 
$\mathcal S_{n,\alpha}^{t_0}$\hspace{.02in}) and there holds that
\begin{equation}\label{ineq1}
0\leq |\sigma(X)-
(t_0\hspace{.02in}\kappa_2+(1-t_0)\kappa_1)|-\Phi_{n,\alpha}(t_0)<\epsilon.
\end{equation}

Since $\lambda_1$ is a real analytic, and thus a $C^3$ curve, we have that
$\lambda_1^{\theta_k}\to\lambda_1^{t_0}$ in $C^{2,1}$. We may then admit  
the existence of  $\nu^{\theta_k}\in \mathcal{S}_{n,\alpha}^{\theta_k}$ such that 
$\lambda_1^{\theta_k}\nu^{\theta_k}\to \lambda_1^{t_0}\mu_j^{t_0}$ in $C^{2,1}$  when $k\to \infty$.

Thus, by continuity, 
for $k$ sufficiently large there are diffeomorphisms 
$X^{\theta_k}\in \Omega^{\lambda_1^{\theta_k}\nu^{\theta_k},\kappa}_{n,\alpha}$ as close as we wish to
$X$ in $C^{2,1}$. Indeed, we may obtain $X^{\theta_k}$ from $X$ by small $C^{2,1}$ perturbations of $X$ in a neighborhood of $\lambda_2$  in order that the trace of $\nu^{\theta_k}$ be the closure of $X^{\theta_k}(\lambda_2)\backslash \lambda_1^{\theta_k}$.  So we may assume that $\tau_*^1(X^{\theta_k})<\kappa$ \hspace{.02in} (since $\tau_*^1(X)<\kappa)$, thus $\nu^\theta_k\in \mathcal S_{n,\alpha}^{\theta_k}$, and have that
\begin{equation}\nonumber
|\sigma(X^{\theta_k})-
(\theta_k\hspace{.02in}\kappa_2+(1-\theta_k)\kappa_1)|-\Phi_{n,\alpha}(t_0)
<\epsilon
\end{equation}
(compare with  ineq. (\ref{ineq1})) when $k$ is sufficiently large.

On the other hand, by (\ref{eep}) and the definition of the function $\Phi_{n,\alpha}$ we have that   the following inequalities must be verified for $k$ large enough:
\begin{equation}\nonumber
|\sigma(X^{\theta_k})-
(\theta_k\hspace{.02in}\kappa_2+(1-\theta_k)\kappa_1)|\geq \Phi_{n,\alpha}(\theta_k)>\Phi_{n,\alpha}(t_0)+\epsilon,
\end{equation}
a contradiction.\bigskip

\noindent
{\bf Case 2:} assume that $\lim_{k\to\infty}\Phi_{n,\alpha}(\theta_k)<\Phi_{n,\alpha}(t_0).$ Let now
\begin{equation}\label{eep1} 
2\epsilon=\Phi_{n,\alpha}(t_0)-\lim_{k\to \infty}\Phi_{n,\alpha}(\theta_k).
\end{equation}
From the definitions of $\xi^t_{n,\alpha}$ and $\Phi_{n,\alpha}$, there holds that for all $k$ there exists $\mu^{\theta_k}\in \mathcal{S}_{n,\alpha}^{\theta_k}$ such that
$0\leq\xi^{\theta_k}_{n,\alpha}(\lambda_1^{\theta_k}\mu^{\theta^k})-\Phi_{n,\alpha}(\theta_k)<
\epsilon/2^{k}$. Hence, from the definition of the function $\xi^t_{n,\alpha}$ there exists
$X^{\theta_k}\in \Omega^{\lambda_1^{\theta_k}\mu^{\theta_k},\kappa}_{n,\alpha}$ 
 such that  
\begin{equation}\nonumber
0\leq |\sigma(X^{\theta_k})-(\theta_k\hspace{.02in}\kappa_2+(1-\theta_k)\hspace{.02in}\kappa_1)|-\Phi_{n,\alpha}(\theta_k)
<\epsilon/2^k
\end{equation} 

By Remark  \ref{rbas} (see also  Remark \ref{resc}), and taking a subsequence if necessary, we may assume that there exists an arc $\mu$ such that the sequence $(X^{\theta_k})$ converge in $C^{2,\alpha}$ to $X\in  \Omega_{n}^{\lambda_1^{t_0}\mu}$  
when $k$ tends to infinity: just take $X$ as the limit in $C^{2,\alpha}$ of some subsequence of $(X^{\theta_k})$ and let the trace of $\mu$ be the closure of $X(\lambda_2)\backslash \lambda_1^{t_0}$.  By Remark \ref{rbas}, we 
also have that $X\in\mathcal X^1$ with $\tau_*^1(X)\leq \kappa$,  because     $\tau_*^1(X^{\theta_k})\leq \kappa$ for any $k\in \mathbb N$. 
Thus we have that
\begin{equation}\nonumber
|\sigma(X)-(t_0\hspace{.02in}\kappa_2+(1-t_0)\hspace{.02in}\kappa_1)|
= \lim_{k\to\infty}\Phi_{n,\alpha}(\theta_k).
\end{equation}

 Then it follows from Lemma $\ref{laux}$  that there exists sequences $(\nu_j)\subset S^{t_0}_{n,\alpha}$ and 
$(X_j)\subset \Omega_{n,\alpha}^{\lambda_1^{t_0}\nu_j}$ such that both $\lambda_1^{t_0}\nu_j\to\lambda_1^{t_0}\mu$ and $X_j\to X$  in $C^{2,1}$. Hence, if $j_0$ is sufficiently large, we have that
\begin{equation}\nonumber
 |\sigma(X_{j_0})-(t_0\hspace{.02in}\kappa_2+(1-t_0)\hspace{.02in}\kappa_1)|< |\sigma(X)-(t_0\hspace{.02in}\kappa_2+(1-t_0)\hspace{.02in}\kappa_1)|+\epsilon.
\end{equation}
Thus,
\begin{equation}\nonumber
\Phi_{n,\alpha}(t_0) \leq\xi_{n,\alpha}^{t_0}(\lambda_1^{t_0} \nu_{j_0} ) \leq |\sigma(X_{j_0})-(t_0\hspace{.02in}\kappa_2+(1-t_0)\hspace{.02in}\kappa_1)| %\label{eqa2}
<  \lim_{k\to\infty}\Phi_{n,\alpha}(\theta_k)+\epsilon.
\end{equation}
A contradiction (see ineq. (\ref{eep1})), and case 2 follows.  \qed
\bigskip

We  finish below the proof that the functions $\Phi_{n,\alpha}$ are  continuous, which is the objective of this section.

\begin{lemma}\label{lcon2} For every  $n\in \mathbb N$ and $\alpha \in (0,1)$ there holds that 
\begin{equation}\nonumber
 \lim_{t\to 0^+}\Phi_{n,\alpha}(t)=\lim_{t\to 1^-}\Phi_{n,\alpha}(t)=0. 
 \end{equation}
\end{lemma}

\proof Take $\nu^t$ as an arc whose trace is the closure of $\lambda_1\backslash \lambda_1^t$. Then $\nu^t\in \mathcal S^t_{n,\alpha}$ for every $n\in \mathbb N$ and $\alpha \in (0,1)$.  We observe that 
\begin{equation}\nonumber
\xi^t_{n,\alpha}(\lambda_1^t\nu^t)\leq |\sigma(X_{12})-(t\hspace{.02in}\kappa_2+(1-t)\hspace{.02in}\kappa_1)|. 
\end{equation}
Hence, 
\begin{equation}\nonumber
0\leq\lim_{t\to 0}\Phi_{n,\alpha}(t)\leq \lim_{t\to 0}\xi^t_{n,\alpha}(\lambda_1^t\nu^t)=0.
\end{equation}

On the other hand, fixed $n\in \mathbb N$ and $\alpha\in (0,1)$, for every $t\in(0,1)$ let $\mu^t\in \mathcal S^t_{n,\alpha}$ be the arc  provided by Lemma \ref{lfal}. Hence, the closed curves $\lambda_1^t\mu^t$ converge to $\lambda_2$ in $C^{2,1}$ (and thus in $C^{2,\alpha}$) when $t\to 1$. When $t$ is close to one we may  obtain $X^t\in \Omega_n^{\lambda_1^t\mu^t}$ by arbitrarily small $C^{2,1}$ perturbations of the identity map $I$ around $\lambda_2$ (this is possible because $\lambda_1^t\mu^t$ converges to $\lambda_2$ in $C^{2,1}$ when $t\to 1$, see Lemma \ref{lfal}). Thus both $X^t$ converges to the identity map $I$ and $X^t(N_2)$ converges to $N_2$ in $C^{2,1}$ when $t$ tends to $1$. Hence, $\lim_{t\to 1} \sigma(X^t)=\sigma(I)=\kappa_2$ and
\begin{equation}\nonumber
0\leq\lim_{t\to 1}\Phi_{n,\alpha}(t)\leq\lim_{t\to 1}|\sigma(X^t(N_2))-(t\hspace{.02in}\kappa_2+(1-t)\hspace{.02in}\kappa_1)|=0.
\end{equation}
\qed

\section  {Proof of the main results}\label{spmr}
\bigskip

We now prove theorem \ref{tp2}. \bigskip

\noindent
{\it Proof of theorem \ref{tp2}.} The hypotheses of Theorem \ref{tp2} except for condition (1) were already assumed by the previous results. We now assume that the first eigenvalue of $N_2$ is equal to $2$ and that the first eigenvalue of $N_1$ is $\kappa_1< 2$.

Fixed $\alpha\in(0,1)$, let $\gamma(t)=\sigma(X_t)$, where $X_t$ belongs to $\Omega_{n,\alpha}^{\lambda_1^t\mu}$ for some $\mu\in \mathcal S_{n,\alpha}^t$ and also satisfies
$\Phi_{n,\alpha}(t)=|\sigma(X_t)-(2\hspace{.03in}(1-t)+t\hspace{.03in} \kappa_1)|$ (for the existence of $X_t$, see the proof of Lemma \ref{lcon1} or the proof of the existence of the diffeomorphisms $X_{n,\gamma}$ below.)
It then holds that $\gamma(t)\leq 2(1-t)+t\hspace{.03in} \kappa_1$: observe that
$\sigma(X_{12}^{-1})<2(1-t)+t\hspace{.03in} \kappa_1$ whenever $t\neq 1$. From Lemma \ref{caux} and Corollary \ref{caux}, there exists a continuous path $\rho:[0,1]\to \mathcal X^1$ such that $\rho(0)=X_t$, 
$\rho(1)=X_{12}^{-1}$, and
$\rho(\theta)\in\Omega_{n,\alpha}^{\lambda_1^t\mu_\theta,\kappa}$ for a compatible arc $\mu_\theta$ and $\theta\in[0,1]$. If $$\sigma(\rho(0))> 2(1-t)+t\hspace{.03in} \kappa_1,$$ then $\Phi_{n,\alpha}(t)>0$ and there exists $\theta\in[0,1]$ such that $\sigma(\rho(\theta))= 2(1-t)+t\hspace{.03in} \kappa_1$, a contradiction (in this case, $\Phi_{n,\alpha}(t)=0$.) Hence, 
$$\gamma(t)=(2\hspace{.03in}(1-t)+t\hspace{.03in} \kappa_1)-\Phi_{n,\alpha}(t).$$ Consequently, from the continuity of the function $\Phi_{n,\alpha}$, which is proved in lemmas \ref{lcon1} and    \ref{lcon2}, the function $t\to\gamma(t)$ is continuous.

Thus,  for every $n\in \mathbb N$ and $\gamma\in(\kappa_1,2)$ there exists 
$t_{n,\gamma}\in(0,1)$ such that 
\begin{equation}
\Phi_{n,\alpha}(t_{n,\gamma})=|\gamma-(2\hspace{.03in}(1-t_{n,\gamma})+t_{n,\gamma}\hspace{.03in} \kappa_1)|.
\end{equation}
Hence, there exists a sequence $\mu^j_{n,\gamma}\in \mathcal S_{n,\alpha}^{t_{n,\gamma}}$ such that $\xi_{n,\alpha}^{t_{n,\gamma}}(\lambda_1^{t_{n,\gamma}}\mu^j_{n,\gamma})$ converges to $\gamma$ when $j$ tends to infinity. 
Denote by $l_{n,\gamma}^j$ the closed curve $\lambda_1^{t_{n,\gamma}}\mu^j_{n,\gamma}$.
Hence there exists a sequence of diffeomorphisms $X_{n,\gamma}^j\in\Omega_{n,\alpha}^{\kappa,\hspace{.03in}l_{n,\gamma}^j}$ such that $\sigma(X_{n,\gamma}^j)$ converges to $\gamma$ when $j$ tends to infinity. By Remark \ref{rbas}, and taking a subsequence if necessary, we may assume that the sequence $(X_{n,\gamma}^j)$ converges in $C^{2,\alpha}$ to a diffeomorphism $X_{n,\gamma}\in \mathcal X^1$ such that $\tau_*^1(X_{n,\gamma})\leq\kappa$. From the continuity of the function $\sigma$ (Lemma \ref{lsigma}) there follows that $\sigma(X_{n,\gamma})=\gamma$.

Fixed $\gamma$ let $(X_{n,\gamma})$ be the sequence of the diffeomorphisms defined above. We may
 again assume by Remark \ref{rbas} and Lemma \ref{lsigma} that a subsequence $(X_{n_j,\gamma})$ converges in $C^{2,\alpha}$ to a diffeomorphism $X_\gamma\in \mathcal X^1$ 
such that $\tau_*^1(X_\gamma)\leq\kappa$ and $\sigma(X_\gamma)=\gamma$.

We now recall the sets of vertices  $V_n\subset N_2$ (see section \ref{s2}), whose union is dense in $N_2$. Since $X_\gamma$ is a diffeomorphism from $N_2$ to $X_\gamma(N_2)$, and diffeomorphisms take dense sets to dense sets then $X_\gamma(\bigcup_{n=1}^\infty V_n)$ is dense in $X_\gamma(N_2)$. Now we recall that $V_{n_j}\subset V_{n_l}$ whenever $j\leq l$. Thus if $p\in \bigcup_{n=1}^\infty V_n$ then there exists $j_0$ such that $p\in V_{n_j}$ for every $j>j_0$. We recall that the mean curvature of $X_{n_j,\gamma}$ at $X_{n_j,\gamma}(p)$ must be equal to zero for every $j>j_0$. Since $X_{n_j,\gamma}$ converges to $X_\gamma$ in $C^{2,\alpha}$ then $X_{n_j,\gamma}(p)$ converges to $X_{\gamma}(p)$ and the mean curvature 
of $X_\gamma$ at  $X_\gamma(p)$ is equal to zero. So $X_\gamma(N_2)$ is a compact embedded $C^{2,\alpha}$ surface with mean curvature equal to zero in a dense set of points and thus it is minimal.

Let then $\gamma_j=2-1/j$. Once more we invoke Remark \ref{rbas} and obtain that a subsequence $X_{\gamma_{j_n}}$ converges in $C^{2,\alpha}$ to $X$, a $C^{2,1}$ diffeomorphism of $A_2$ such that $\tau_*^1(X)\leq\kappa$. Let $M_n= X_{\gamma_{j_n}}(N_2)$. The sequence of closed embedded minimal surfaces $M_n$ then converges in $C^{2,\alpha}$ to a closed embedded minimal surface $M=X(N_2)$.  Moreover, the first eigenvalue of $M_n$ is equal to $2-1/j_n$ (so these surfaces are non-congruent) and, by Lemma \ref{lsigma}, the first eigenvalue of $M$ is two.\qed
\medskip

Theorem \ref{tp1} is proved by a reasoning that  is fully described in \cite{p1} so we will just point out its main features. Indeed, the argument in \cite[Section 4]{p1} does not require that the closed surfaces that are dealt with are tori and its extension to surfaces of arbitrary genus is straightforward. Proofs of theorems \ref{tp1} and \ref{tp2} are then similar in structure and contents, their main difference being  the definition of the function $\xi$, which is adapted to the specific result that is proved. 

Finally, we want to stress that Theorem \ref{tp1} is not weaker than Theorem \ref{tp2}. Although Theorem \ref{tp1} does not control the eigenvalues of the sequence of closed embedded  minimal surfaces whose existence it guarantees, it assures that this sequence does converges to $N_2$ in $C^{2,\alpha}$. Thus $N_2$ is not isolated in the set of closed embedded minimal surfaces in $\S3$. On the other hand, Theorem \ref{tp2} does not permit such conclusion. 

\section{Final remarks}\label{sfr}

Theorems \ref{tp1} and \ref{tp2} are meant to be used in settling some issues of the space of embedded minimal surfaces in $\S3$. 

An application that we devise for Theorem \ref{tp1}  is a proof that either there exists, up to congruences, only one minimal surface of a fixed genus $g$ embedded in $\S3$ or any of such surfaces is not isolated in $C^{2,\alpha}$, i.e., given $M$, a closed minimal surface embedded in $\S3$ of genus $g$ there exists a sequence $M_n$ of non-congruent closed minimal surfaces embedded in $\S3$ converging to $M$ in 
$C^{2,\alpha}$.       

As for Theorem \ref{tp2}, as we have already put forward in the introduction to this work, it is conceived to be used to prove the well known Yau's conjecture (or answer the corresponding question in the affirmative) that every minimal closed surface embedded in $\S3$ is embedded by the first eigenvalue.

In order to prove these results, we must first verify that conditions 2 and 3 of theorems \ref{tp1} and \ref{tp2} are satisfied by any possible pair $N_1$, $N_2$ of non-congruent closed minimal surfaces embedded in $\S3$ of genus $g$. Thus, it must be proved that there exists embedded closed curves $\lambda_1\subset N_1$ and $\lambda_2\subset N_2$ contained in equators of $\S3$ satisfying such requisites. 

A possible way to verify the existence of these curves is to apply and generalize Section 3 of \cite{p1} to embedded minimal surfaces of arbitrary genus. This consists in using the {\it two-piece} property for embedded minimal surfaces discovered by Ros (see \cite{rs}), i.e., that every equator divides a closed surface embedded in $\S3$ in exactly components, to prove the existence of the curves $\lambda_i$ with the prescribed properties as is done in \cite{p1} for tori.

Satisfied its conditions, Theorem \ref{tp1} naturally leads to its application suggested above. On the other hand, a proof of Yau's conjecture requires an additional  step: we  have to prove that any closed minimal surface embedded by the first eigenvalue is isolated in the $C^{2,\alpha}$ topology from embedded closed minimal surfaces that are not embedded by the first eigenvalue. This is shown in Section 5 of \cite{p1} for tori by a reasoning whose generalization to surfaces of arbitrary genus is immediate. 
Hence, the sequence built in Theorem \ref{tp2} can not exist and we reach a contradiction: satisfied the necessary conditions, $N_1$ must be embedded by the first eigenvalue as well.

Another way to prove the results above may be dropping from the theorems $\ref{tp1}$ and $\ref{tp2}$ the condition on the existence of curves $\lambda_1$ and $\lambda_2$ as above. The control on the deformation $t\to X_t$ from $X_{12}^{-1}$ to the identity map would use instead the distance function $\tau^\alpha_*$, or the absolute value 
$|\sigma(X_t)-(2(1-t)+t\kappa_1)|$, or a similar function.

\vspace{.35in}

\centerline{\large \bf APPENDICES} 

\appendix

\vspace{.1in}

\section{Statement and proof of Lemma \ref{laux}}\label{A1}

\vspace{.2in}

We conform here to the definitions and notations established in Section \ref{s1}.
We then assume that $\mathcal X^\alpha/\sim$, $0<\alpha\leq 1$, is endowed with the metrics $d_*^\alpha$ and let  $\tau_*^\alpha(X)=d_*^\alpha(\bar X,\bar I)$ for $X\in \mathcal X^\alpha$.

We remark that $\mathcal X^\alpha$ is a Banach manifold, indeed a submanifold of a Banach space. 
In order to verify this, let $V$ be the space of $C^{2,\alpha}$ maps $Z:\S3\to\mathbb R^4$, $0<\alpha\leq 1$, with the holder $C^{2,\alpha}$ norm described in Section \ref{s1} by eqs. (\ref{eqnh1}) and (\ref{eqnh2}). The vector space $V$ is then a Banach space. Let $M_1\subset V$ be the set of diffeomorphisms of $\S3$ into their images. Hence, $M_1$ is an open subset (and thus a submanifold) of $V$. Let $M_2\subset V$ be the set of applications of $\S3$ into itself. Thus $M_2$ is also  a submanifold of $V$. We then have that $M=M_1\cap M_2$, i.e., the set of $C^{2,\alpha}$ diffeomorphisms of $\S3$, is a submanifold of $V$.

Let $U=\{(Z,Z^{-1})\hspace{.02in}|\hspace{.02in}Z\in M\}$. The space $U$ is then diffeomorphic to $M\times M$, thus the set $U$ is a submanifold of the product $V\times V$, a Banach space with the product metric.  Finally, we observe that
the space $\mathcal X^\alpha$, whose metric is given by Definition \ref{ddist}, is identical to $U$, and is thus a submanifold of a Banach space.  

We shall now obtain some submanifolds of $\mathcal X^\alpha$ itself.

Firstly, fix a unitary vector field $n$ normal to $N_2$ in $\S3$ and let $X\in \mathcal X^\alpha$. We may then orient $X(N_2)$ by $dX(n)$ and thus define a mean curvature function along $X(N_2)$.
If $q$ is a point in $N_2$ we let $H_q(X)$ be the mean curvature in $\S3$ of $X(N_2)$ at $X(q)$. Hence, the function $H_q: \mathcal X^\alpha\to \mathbb R$ is differentiable and $H_q^{-1}(0)$ is a submanifold of $X^\alpha$ because {\it zero} is a regular value of $H_q$.

We then remark that the set $W_t^\alpha$ of diffeomorphisms $X\in \mathcal X^\alpha$ such that $X(\lambda_2)$ is in an equator of $\S3$ and  contain a subset congruent to $\lambda_1^t$ is also a submanifold of $X^\alpha$.

Finally, let $\mathcal S^t$ be the set of closed arcs $\mu$ such that $\lambda_1^t\mu$ is a simple $C^{2,1}$ curve in a equator. We recall the definition of $V_n$ (see Definition \ref{dvn}).
Then the set $\Lambda_n^t=\bigcup_{\mu\in\mathcal S^t}\Omega_n^{\lambda_1^t\mu}$   is a submanifold of $\mathcal X^1$, since it is the intersection $\left(\hspace{.02in}\bigcap_{q\in V_n}\hspace{.02in}
 H_{q}^{-1}(0)\hspace{.01in}\right)\cap W_t^1$ of the transversal submanifolds defined above.

\begin{lemma}\label{laux} Let $X\in\Omega_n^{\lambda_1^t\mu}$ for some $\mu\in \mathcal S^t$ and $\tau_*^1(X)=\kappa$. For every $j$ there exists   $\mu_j\in S^t_{n,\alpha}$ and $X_j\subset\Omega_{n,\alpha}^{\lambda_1^t\mu_j,\kappa}$ such that the sequence $(X_j)$ converges to $X$ in $C^{2,\alpha}$.
\end{lemma} 

\begin{remark} If $Z\in \mathcal X^\alpha$ then one can perturb $Z$ and $Z^{-1}$ in the points where the functions that compose the norms $||Z-I||_{C^{2,\alpha}}$ and $||Z^{-1}-I||_{C^{2,\alpha}}$ take their maximal values. So if  $Z^*$ is the  perturbation of $Z$ in $\mathcal X^\alpha$ and the perturbation is small enough then we may assume that 
\begin{equation}\label{eaf}
|\tau^1_*(Z^*)-\tau^1_*(Z)|>\delta\hspace{.02in} d^\alpha_*(Z^*,Z) 
\end{equation}
for a given constant $0<\delta<1$. 

We may also assume that we can continuously extend this perturbation to a neighborhood $A$ of $Z$ in $\mathcal X^\alpha$ so a diffeomorphism $Y\in A$ is perturbed to a diffeomorphism $Y_*$ such that $|\tau^1_*(Y^*)-\tau^1_*(Y)|>\delta/2\hspace{.03in} d^\alpha_*(Y^*,Y)$.

We may proceed similarly with a diffeomorphism $Z$ in a submanifold $\Lambda_n^t$, perturbing it in $\Lambda_n^t$
to a diffeomorphism $Z^*$ such that $Z$ and $Z^*$ satisfy the inequality (\ref{eaf}), and then extend the perturbation as above to a neighborhood of $Z$ if we assume that $\tau^1_*(Z)$ is sufficiently greater than
$\tau^1_*(X_{12})$. In this case, $\tau^1_*(Z)$ shall be  larger than the curvatures of $\lambda_1$, $\lambda_2$ and their derivatives, thus the perturbation above would not have to  alter the curvature of $\lambda_1^t$ in order to have ineq. (\ref{eaf}) satisfied and the perturbation extended as above.

Since the only restriction on $\kappa$ is being larger than $\tau^1_*(X_{12})$,  we may assume the existence of the perturbation above for the diffeomorphism $X$ of Lemma \ref{laux}.
\end{remark}

\noindent
{\bf Proof of Lemma \ref{laux}.} Let $A$ be a small neighborhood of $X$ in the Banach manifold $\Lambda_n^t$ (defined above). We can then identify $A$ to an open subset of a  Banach space $\bf E$. Let $V$ be a unitary vector in $\bf E$ and $\epsilon$, $\delta>0$ 
such that if $0<\theta<\epsilon$ then $$\tau^1(X+\theta V)>\delta\hspace{.04in}(1+\theta)\hspace{.02in}\tau^1(X).$$
We may also assume that $V$ is chosen in such a way that if $A$ is small enough, $Y\in A$, and $0<\theta<\epsilon$ then $\tau^1(Y+\theta V)>\delta/2\hspace{.02in}(1+\theta)\hspace{.02in}\tau^1(Y)$. Hence, by the Generalized Picard-Lindelof Theorem (see \cite[p. 78, Theorem 3.A]{zd}), the initial value problem $x(0)=X,$ $x'(t)=V$, has a continuous differentiable solution in a interval $(-a,a)$ for some $0<a<\epsilon$.

We then observe that if $t<0$, then $\tau^1(x(t))<\tau^1(X)=\kappa$. Moreover, when $t$ tends to $zero$, $x(t)$ tends to $X$ in $C^{2,1}$, thus in $C^{2,\alpha}$.  \qed

\begin{corollary}\label{caux} The set $\{X\in \Lambda_n^t\hspace{.03in}:\hspace{.03in}\tau^1_*(X)\leq \kappa\}$
is star-shaped.
\end{corollary}

\vspace{.3in}

\section{Statement and proof of Lemma \ref{lfal}}\label{A2}

Let $\varphi: \mathbb R\to \mathbb R$ be a non-negative $C^\infty$ function such that $\varphi(x)=0$ if $|x|\geq 1$ and $\varphi(x)=1$ if $|x|\leq1/2$. Let then $\varphi_{\epsilon}(x)=\varphi(x/\epsilon)$.
For $k=0,1,\dots$, let $M_k=\max|\varphi^{(k)}(x)|$.

\begin{lemma}\label{lfal0} Let $f:\mathbb R \to \mathbb R$ be an analytic function such that $f(0)=f'(0)=f''(0)=f'''(0)=0$. Hence, $h_\epsilon(x)=\varphi_\epsilon(x)f(x)$ converges in $C^3$ to the identically null function when $\epsilon\to 0$.
\end{lemma}

\proof  The hypotheses above imply that the function $f$ has the ensuing expansion around zero: $$f(x)=f^{(4)}(0)\hspace{.02in}\frac{x^4}{4!}+O(x^5).$$ 
We also have that $h_\epsilon(x)= 0$ when $|x|\geq \epsilon$. If $|x|<\epsilon$ then the following inequalities hold:

$$\left|{h_\epsilon'}{(x)}\right|= \left|\hspace{.02in}\hspace{.02in}\varphi(\frac{x}{\epsilon})\left(f^{(4)}(0)\frac{x^3}{3!}+O(x^4)\right)+ \frac{1}{\epsilon}\hspace{.03in}\varphi'(\frac{x}{\epsilon})\left(f^{(4)}(0)\frac{x^4}{4!}+O(x^5)\right)\hspace{.01in} \right|\leq$$
$$ M_0\left( |f^{4}(0)|\frac{\epsilon^3}{3!}+O(\epsilon^4)\right)+M_1\left( |f^{(4)}(0)|\frac{\epsilon^3}{3!}+O(\epsilon^4)\right)\leq O(\epsilon^3).$$

In the same way we obtain that  $|h_\epsilon''(x)|\leq O(\epsilon^2)$ and $|h_\epsilon'''(x)|\leq O(\epsilon)$, thus proving  Lemma \ref{lfal0}. \qed

\begin{lemma}\label{lfal} For any $t\in(0,1)$ we can choose an arc $\mu^t\in S_{n,\alpha}^t$ such that  $\lambda_1^t\mu^t$ converges to $\lambda_2$ in $C^{2,1}$ (and thus in $C^{2,\alpha}$) when $t\to 1$.
\end{lemma}

\proof Along this work we identify subsets of $\S3$ up to congruences. We can then assume that $\lambda_1$ and $\lambda_2$ are in a same equator $S(v)$ of $\S3$, that $\lambda_1(0)=\lambda_2(0)$, and that $\lambda_1'(0)=\lambda_2'(0)$, where we take both curves parameterized by arc length in $\S3$.
Let $p\in S(v)$ be a point not in $\lambda_1$ or $\lambda_2$, and $x: S(v)\backslash\{p\}\to\mathbb R^2$ be an analytic parameterization of $S(v)\backslash \{p\}$ (which may be obtained, e.g., by stereographic projection). Thus we may identify $\lambda_1(\theta)=(f_1(\theta),g_1(\theta))$ and $\lambda_2(\theta)=(f_2(\theta),g_2(\theta))$, which are real analytic curves in $\mathbb R^2$. We then let
$$\lambda_1^t\mu^t(\theta)= \lambda_2(\theta)+h_{2t}(\theta)\left(\lambda_1(\theta)-\lambda_2(\theta)\right),$$
for $\theta \in (-L_2/2,L_2/2)$, where $h_{\epsilon}$ is the  function defined in Lemma \ref{lfal0} and $L_2$ is the length of $\lambda_2$. 
We remark we obtain  $\mu^t$ from its trace, i.e., the closure of $\lambda_1^t\mu^t ((-L_2,L_2))\backslash \lambda_1^t\mu^t ((-t,t))$. We also observe that $\lambda_1^t\mu^t$ is a $C^\infty$ embedded closed curve when $t>0$ is sufficiently small.

We now consider the difference: 
$$\lambda_2(\theta)-\lambda_1^t\mu^t(\theta)=h_{2t}(\theta)\left(\lambda_2(\theta)-\lambda_1(\theta)\right).$$
We already have that $\lambda_1(0)=\lambda_2(0)$ and $\lambda_1'(0)=\lambda_2'(0)$. We remind that  a third order contact between two curves is preserved by $C^3$ diffeomorphisms. The conditions on the curvatures of $\lambda_1$ and $\lambda_2$ at $t=0$ (at $t=0$, these curvatures are equal and  assume critical values for their respective curves, see Section \ref{sst}) then imply that $\lambda_1''(0)=\lambda_2''(0)$ and
$\lambda_1'''(0)=\lambda_2'''(0)$ as well. Hence, by Lemma \ref{lfal0}, the function $\theta \to h_{2t}(\theta)\left(\lambda_2(\theta)-\lambda_1(\theta)\right)$ converges uniformly up to third order to the identically null function when $t$ goes to $zero$, which implies Lemma \ref{lfal}.\qed 
\medskip\bigskip

\noindent
Fernando A. A. Pimentel ({\small pimentelf@gmail.com})

\end{document}